\documentclass[12pt]{article}
\usepackage{amsfonts,amsmath,amscd,amsthm}

\newtheorem{prop}{Proposition}[section]
\newtheorem{lem}[prop]{Lemma}
\newtheorem{thm}[prop]{Theorem}

\newtheorem{exe}[prop]{Example}

\font\bulo= eusm9 at 12pt
\font\bbulo= cmssbx10 at 12pt

\newcommand{\AV}{ \mbox{\bulo A} }

\newcommand{\AVA}{ \mbox{\bulo F} }
\newcommand{\AVAP}{ \AVA\,' }
\newcommand{\SH}{ \mbox{\bbulo H} }
\newcommand{\EM}{ \mbox{\bulo E} }
\newcommand{\IS}{ \mbox{\bulo I} }
\newcommand{\ISP}{ \IS\,' }
\newcommand{\MOR}{ \mbox{\bulo F} }
\newcommand{\MORP}{ \MOR\,' }
\newcommand{\tp}{ {^t}\!}

\newcommand{\X}{ \mbox{\bulo X} }
\newcommand{\Y}{ \mbox{\bulo Y} }
\newcommand{\V}{ \mbox{\bulo V} }
\newcommand{\W}{ \mbox{\bulo W} }

\newcommand{\C}{\mathbb C}
\newcommand{\R}{\mathbb R}
\newcommand{\Z}{\mathbb Z}
\newcommand{\Q}{\mathbb Q}

\newcommand{\Ker}{ \mbox{Ker} }

\begin{document}

\title{Siegel coordinates and moduli spaces \\ for 
morphisms of Abelian varieties}
\author{Lucio Guerra}
\date{}
\maketitle

\begin{abstract}
\noindent We describe the moduli spaces of morphisms
between polarized complex abelian varieties. The discrete
invariants,  derived from a Poincar\'e decomposition of
morphisms,  are the types of polarizations and of
lattice homomorphisms occurring in the decomposition.
For a given type of morphisms the moduli variety is irreducible, and is obtained 
from a product of Siegel
spaces modulo the action of a discrete group. 

\end{abstract}

\section{Introduction}

We describe the moduli spaces of morphisms
between polarized abelian varieties over the complex field.
Several instances of these modular varieties occur in the literature.
The moduli spaces of isogenies between elliptic curves occur as Hecke correspondences on modular
curves, cf. \cite{KM}, \cite{Mi}. The moduli spaces of embeddings of elliptic
curves into abelian surfaces  appear through their  images in the
moduli spaces of abelian surfaces. The moduli spaces of abelian surfaces containing elliptic curves
coincide with certain Humbert modular surfaces, as shown by Kani
\cite{K}, and may also be seen as spaces of curves of genus two which cover elliptic curves, as in
Murabayashi \cite{Mu}.  Moduli spaces of isogenies on abelian surfaces have been studied by Hulek and
Weintraub \cite{HW}, and by Birkenhake and Lange \cite{BL}.

In this paper we propose a general treatment, 
in terms of Siegel coordinates. We identify  the discrete invariants which
correspond to the irreducible components of the moduli space of morphisms, and we describe the
irreducible components as quotients of products of Siegel spaces. 

The discrete invariants 
are derived from a {Poincar\'e decomposition of
morphisms}. A morphism
$V\rightarrow W$ of polarized abelian varieties, 
after suitable isogenies
$X\times X'\rightarrow V$ and $Y\times
Y'\rightarrow W$, each a sum of two embeddings
of complementary abelian subvarieties, lifts to a morphism
$X\times X'\rightarrow Y\times Y'$ which is
given by an isogeny
$X\rightarrow Y$ multiplied by the zero map on $X'$.

The sequence $\delta$ which collects
the {types of polarizations} of the six
varieties occurring in the decomposition is an
invariant of the morphism. A second invariant
is defined by taking into account how the
morphisms in the decomposition behave with
respect to the polarizations.

A  morphism of abelian varieties $X\rightarrow
Y$ with respect to symplectic bases gives rise
to an integer matrix $M$ which represents the
induced homomorphism 
$H_1(X,\Z) \rightarrow H_1(Y,\Z)$, and what is only dependent on the
polarizations is the equivalence class 
$[M]$ under the natural action of
symplectic groups. In the same way 
a Poincar\'e decomposition, which consists of a sequence
of three morphisms, determines a
sequence $\tau$ of integer matrices relative
to symplectic bases, and
what is invariant is an equivalence class $[\tau]$, 
defined in the natural way.
The collective datum \
$\delta,[\tau]$ \ is what we call a
{\em type of morphisms} between polarized
abelian varieties, and we characterize 
the discrete data which are types of
morphisms.   

The global coordinates for morphisms  are  the Siegel coordinates of the
varieties occurring in the Poincar\'e decomposition. As the type is fixed, it turns out that the
morphism is determined if the varieties $X,X'$ and $Y'$ are given. Thus,
for every type of morphisms 
we construct a complex analytic variety \
$\AVA_{\delta}[\tau]$ \ which is a coarse moduli
space for morphisms of the given type. 
The moduli variety is irreducible, 
and is obtained from a product of Siegel
spaces modulo the action of a discrete group.

The construction is presented in several steps.
In \S \ref{firstinvariant}
we introduce the invariant $[M]$ 
mentioned before. Using this, in \S \ref{isogenies} the moduli
space of isogenies is constructed. In \S \ref{Hecke} we describe some known examples,
such as the Hecke correspondences. 
In \S \ref{refinedinvariant} we introduce a refined invariant, the type of an
embedding, derived from classical Poincar\'e
reducibility. In \S \ref{embeddings} the moduli space of
embeddings is constructed. In \S \ref{reducibility} the Poincar\'e
reducibility for morphisms is presented and in \S \ref{generalmorphisms}
the moduli space of general morphisms is
constructed. Finally in \S \ref{proofs} we give the proof
of the moduli property for these spaces, in a
sketched form, and in \S \ref{questi} we shortly discuss the question of classifying the discrete
data which are types of morphisms.

\newpage

\section{A first invariant} \label{firstinvariant}

Some notation.
Every abelian variety $X$ will be endowed with a polarization $L_X$, that we identify with the
corresponding alternating form on $H_1(X,\Z)$. We denote by $D = (d_1,\ldots,d_n)$ the
sequence of elementary divisors, the type of the polarization, where $n$ is the dimension of $X$.
The alternating form will be determined by means of some symplectic basis $\lambda = (\lambda_1,
\ldots, \lambda_{2n})$, such that $L_X (\lambda_i, \lambda_j)$ is equal to  $d_i$ for $j=i+n$, to
$-d_i$ for $i=j+n$, and 0 otherwise,  two symplectic bases for the polarization being related by some
element of the symplectic group $Sp(D,\Z)$. We use the symbol
$\widehat{D}$ for the matrix of the alternating form with respect
to a symplectic basis $\lambda$, and  sometimes we write this as $\langle\lambda,\lambda\rangle =
\widehat{D}$, dropping the symbol $L_X$. 
Note that the datum $D$ encodes the dimension $n$.

The collection of isomorphism classes of abelian varieties endowed with  symplectic
basis for a polarization of type $D$ will be denoted by $\AV'_D$. The superscript is to remind the
choice of symplectic bases, and in this way it will be used again in the following. The
bijective map  $\AV'_D \leftrightarrow \SH_n$ with the Siegel  space introduces the structure of
a complex manifold. The collection of isomorphism classes of abelian varieties with polarization of
type $D$ is bijective to the quotient variety $\SH_n / Sp(D,\Z)$. In a more precise language
it is well known that one should speak of varieties representing the appropriate moduli functors.

\bigskip
We consider morphisms $f: X\rightarrow Y$ between abelian varieties, both endowed with polarization.
As morphisms between abelian varieties we mean in the following the morphisms of analytic groups,
the same as the morphisms of analytic varieties up to translations in the range. 

An isomorphism $f \cong f'$ with another morphism $f': X'\rightarrow Y'$ is a couple of isomorphisms
$u: X\rightarrow X'$ and $v: Y\rightarrow Y'$, which preserve the polarizations, such that $f' = v
\circ f \circ u^{-1}$. We denote by $D,E$ the pair of types of polarizations on $X,Y$, 
of dimensions $n,m$. 
\bigskip

\noindent {\bf Definition.}
Define $\MORP_{D,E}$ to be the collection of morphisms $f$ of abelian varieties, both endowed
with polarizations of the given types and with symplectic bases for the polarizations, 
modulo isomorphisms $f \cong f'$ which preserve (the polarizations and) the bases.
Define ${\MOR}_{D,E}$ to be the collection of morphisms of abelian varieties with polarizations of
the given types, modulo isomorphisms (which preserve the polarizations).
Define $\Gamma_{D,E} = Sp(D,\Z)\times Sp(E,\Z)$.
\bigskip

The group $\Gamma_{D,E}$ acts on the set $\MORP_{D,E}$ 
by producing the change of symplectic bases. The element $(A,B)$  acts sending the
isomorphism class represented by the morphism $f$ endowed with symplectic bases $\gamma,\lambda$
into the isomorphism class represented by $f$ with the bases 
$\gamma A$ and $\lambda B$.  There is an identification
$${\MOR}_{D,E} \ = \ \begin{array}{c}\MORP_{D,E} \\ \hline \Gamma_{D,E}
\end{array}$$

To a morphism $f$ is associated the homomorphism
$f_*:H_1(X,\Z)\rightarrow H_1(Y,\Z)$, called the {rational representation} of $f$. With respect
to the symplectic bases $\gamma, \lambda$, it is represented by an integer matrix $M$, that we also
consider as a homomorphism $\Z^{2n} \stackrel{M}{\rightarrow} \Z^{2m}$. This defines a
map
$$\MORP_{D,E} \longrightarrow Hom(\Z^{2n},\Z^{2m})$$
sending the isomorphism class of $(f,\gamma,\lambda)$ to the matrix $M$. We denote by
$\MORP_{D,E}(M)$ the fibre over the matrix $M$.

The group $\Gamma_{D,E}$ also acts on the set $Hom(\Z^{2n},\Z^{2m})$ 
by producing the equivalence of matrices that represent a given homomorphism $f_*$
under a change of bases. The element $(A,B)$ acts as $M \mapsto BMA^{-1}$. 
The map above is equivariant with respect to the actions of $\Gamma_{D,E}$, 
so there are isomorphisms   
of the fibres over two matrices in the same orbit. 
Consider the quotient map
$${\MOR}_{D,E} = \begin{array}{c} \MORP_{D,E} \\ \hline \Gamma_{D,E} \end{array}
 \ \longrightarrow \ 
\begin{array}{c} Hom(\Z^{2n},\Z^{2m}) \\ \hline \Gamma_{D,E} \end{array}$$
sending the isomorphism class of $f$ to the equivalence class $[M]$.
We denote by  ${\MOR}_{D,E}\,[M]$ the fibre over the equivalence class $[M]$. 
There is an identification
$${\MOR}_{D,E}\,[M] \ = \ \begin{array}{c} \MORP_{D,E}(M) \\ \hline \Gamma_{D,E}(M) \end{array}$$
where
$\Gamma_{D,E}(M)$ denotes the stabilizer of $M$

The equivalence class $[M]$ is a first invariant that may be used in order to
separate the components of the space ${\MOR}_{D,E}$.

\section{Moduli of isogenies} \label{isogenies}

The invariant defined in the preceding section is indeed sufficient 
for classifying the isogenies. 

Let $f:X\rightarrow Y$ be a morphism of abelian varieties of the same dimension $n$, endowed with
polarizations of types $D,E$, and with symplectic bases, and let $M$ be the matrix of $f_*$ relative
to the bases. The morphism is an isogeny if and only if  $\det M\neq 0$. 
Then there is an exact sequence 
\begin{equation} \label{seq}
0\longrightarrow \Z^{2n} \stackrel{M}{\longrightarrow} \Z^{2n} {\longrightarrow} F
\longrightarrow 0
\end{equation}
where the group $F$ is defined as the cokernel of $M$. 
Moreover there is an isomorphism $F\cong\Ker(f)$.

We say that $f$ is an {\em isogeny of polarized abelian varieties} if  $L_X = f^* L_Y$
holds, and this happens if and only if 
\begin{equation} \label{tipi} \widehat{D}  =  \tp M \widehat{E}M \end{equation}
Note that $\det M \neq 0$ follows from this. For instance one has polarizations of the same type if
and only if $M$ is symplectic, and then $F=0$.

We call the datum $\delta = (D,E)$  the {\em polarization type} of the isogeny $f$. We call the datum
$\delta, M$ the {\em type of the isogeny relative to the symplectic bases}. 
\bigskip

\noindent {\bf Definition.}
Let $\delta = (D,E)$ be a pair of types of polarizations, of the same dimension $n$.
Define the subsets $\ISP_\delta \subset \MORP_{D,E}$ and $\IS_\delta \subset \MOR_{D,E}$ 
consisting of isomorphism classes of isogenies of  polarized abelian varieties,
with symplectic bases selected or not selected. Write moreover $\Gamma_\delta = \Gamma_{D,E}$.
\bigskip

So we have the restricted map
$$\ISP_\delta \longrightarrow GL(2n,\Z)$$
The image is  characterized by equality (\ref{tipi}), and
the fibres are described as $\ISP_\delta(M) = \MORP_{D,E}(M)$. 
Then we  have the quotient map
$$\IS_\delta = \frac{\ISP_\delta}{\Gamma_\delta} \ \longrightarrow \ \frac{GL(2n,\Z)}
{\Gamma_\delta}$$ 
sending the isomorphism class of the isogeny $f$ to the equivalence class $[M]$.

We say that the datum
$\delta, [M]$ is the {\em type of the isogeny}  (tout-court, relative to the polarizations,
independent  of symplectic bases).

The fibre of the quotient map  is described as 
$$\IS_\delta[M] \ = \ \dfrac{\ISP_\delta(M)}{\Gamma_\delta(M)}$$
where $\Gamma_\delta(M)$ denotes the stabilizer of $M$.

For a fixed type of isogenies, consider the two
projections
$$\begin{array}{ccc}
& {\ISP_\delta}(M) \vspace{5pt} \\ & \swarrow \hspace{10pt} \searrow \\ 
\AV'_D &  & \AV'_E 
\end{array}$$      
in which to $f:X\rightarrow Y$ are associated the varieties $X$ and $Y$ respectively. 

\begin{prop} \label{duebiiettive}
The two projections are bijective.
\end{prop}

\begin{proof}
There is a bijective correspondence on the level of comples tori, which restricts to abelian
varieties. Given a complex torus $Y$ and a basis in $H_1(Y)$, i.e. an isomorphism $H_1(Y)
\cong \Z^{2n}$, from the  exact sequence (\ref{seq})  a representation
$H_1(Y)\rightarrow F\rightarrow 0$ is deduced, and this determines a complex torus $X$ together with
a covering map $X\rightarrow Y$, unique up to isomorphisms.
Then from the exact sequence an isomorphism $H_1(X) \cong \Z^{2n}$ is also obtained.
Conversely, given a complex torus $X$ and an isomorphism $H_1(X) \cong \Z^{2n}$, 
by means of the exact sequence  an inclusion $F\subset X$ is deduced, 
hence a quotient torus $Y=X/F$, together with an isomorphism
$H_1(Y) \cong \Z^{2n}$.

Let $f:X\rightarrow Y$ be a morphism of complex tori constructed as above.
The basis in $H_1(Y)$ is symplectic for a unique antisymmetric form $L_Y$ of type $E$. 
This determines a natural $\R$-bilinear form $\bar L_Y$ 
on the space ${H^0(\Omega_Y)}^*$ such that $\mbox{im}(\bar L_Y)$ 
coincides with $L_Y$ on $H_1(Y)$, and one has a polarization if and only if $\bar L_Y$
is a positive hermitian form. 
Similarly the basis in $H_1(X)$ is symplectic 
for an antisymmetric form $L_X$ of type $D$, which determines a form $\bar L_X$.
The two bases are related through the matrix $M$.
The relation (\ref{tipi}) between the types is the condition so that there is equality
$L_X = f^* L_Y$ of the antisymmetric forms. This implies $\bar L_X = f^*\bar L_Y$. 
As the induced isomorphism $f^*:{H^0(\Omega_X)}^*\rightarrow{H^0(\Omega_Y)}^*$
is indeed $\C$-linear, it follows that $\bar L_X$ is a 
positive hermitian form if and only if the same holds for $\bar L_Y$.
\end{proof}

The sets $\AV'_D,\AV'_E$ both have the structure of an analytic variety, 
represented by the Siegel space $\SH_n$. Because of \ref{duebiiettive}
there is a diagram of bijections
$$\begin{array}{ccc}
& {\ISP_\delta}(M) & \medskip \\
 & \swarrow\hspace{15pt}\searrow &  \smallskip \\
\SH_n\!\! & \longleftrightarrow &\!\! \SH_n
\end{array}$$

\begin{prop} \label{locclos}
The bijection in the bottom line is an automorphism of the Siegel space.
\end{prop}

\begin{proof}
The Siegel space ${\SH}_n$ may be viewed as a quotient: the space $L(D)$ of lattice bases $\lambda =
(\lambda_1, \ldots, \lambda_{2n})$ in  $\C^n$, subject to the locally closed condition that
the alternating form defined by $\langle\lambda, \lambda\rangle =\widehat D$ is the imaginary part of
some positive hermitian form (a condition that the Riemann relations express in terms of the Siegel
coordinates of $\lambda$), divided by the natural action of $GL(n,\C)$.
An equivariant isomorphism $L(D) \leftrightarrow L(E)$  induces an automorphism of ${\SH}_n$.

If the lattice basis $\lambda$ defines a torus $Y$ then the lattice basis $\gamma:=\lambda M$
determines the complex torus $X$ that covers $Y$; conversely if 
$\gamma$ defines a torus $X$ then $\lambda := \gamma M^{-1}$ defines the quotient torus $Y$.
This correspondence $\gamma \leftrightarrow \lambda$ is an automorphism of the space of all lattice
bases. It follows from the preceding proof that the correspondence restricts to an
isomorphism $L(D) \leftrightarrow L(E)$.
\end{proof}

The stabilizer subgroup ${\Gamma}_\delta(M)$ consists of the pairs  $(A,B) \in \Gamma_\delta$ 
such that $MA=BM$.  In other words there is a diagram of group homomorphisms
$$\begin{array}{ccc}
& {\Gamma}_\delta(M) \vspace{5pt} \\ & ^\alpha\swarrow \hspace{25pt} \searrow\,^\beta \\ 
GL(2n,\Q) & \begin{array}{c}{\longleftarrow\!\longrightarrow}\vspace{-10pt}\\ _{\phi_M}\end{array}
& GL(2n,\Q)
\end{array}$$
where $\alpha,\beta$ are the two projections, and where $\phi_M$ 
is the inner automorphism 
such that $A\mapsto MAM^{-1}=B$. 
Equality (\ref{tipi}) implies that $\phi_M$ restricts to an isomorphism $Sp(D,\Q)
\leftrightarrow Sp(E,\Q)$ of the rational symplectic groups.
So there is an induced isomorphism
$$\begin{array}{ccc}
Sp(D,\Z)\,\cap\,\phi_M^{-1} GL(2n,\Z)  
& \longleftrightarrow &  \phi_M GL(2n,\Z)\,\cap\,Sp(E,\Z)  \\ 
_{||} & & _{||} \vspace{5pt} \\
\alpha\,{\Gamma}_\delta(M) &  & \beta\,{\Gamma}_\delta(M)
\end{array}$$

The actions of the groups $\alpha\,{\Gamma}_\delta(M)$ and $\beta\,{\Gamma}_\delta(M)$ 
on $\SH_n$ are properly discontinuous, because they are restrictions
of the actions of symplectic groups. It is well known that on the quotient spaces there are
natural structures of analytic varieties. Moreover in the diagram of bijections
$$\begin{array}{ccc} & \IS_\delta[M] & \medskip \\ 
& \swarrow\hspace{10pt}\searrow & \\
\begin{array}{c} \SH_n \\ \hline \alpha\,{\Gamma}_\delta(M) \end{array} 
& \longleftrightarrow & 
\begin{array}{c} \SH_n \\ \hline \beta\,{\Gamma}_\delta(M) \end{array}
\end{array}$$
the arrow in the bottom line is an isomorphism of varieties. This determines
on $\IS_\delta[M]$ a unique structure of a complex analytic variety. Being a quotient of $\SH_n$,
this is an irreducible variety.

\begin{thm}\label{isogeniale}
The variety $\IS_\delta[M]$ is a coarse moduli space for isogenies of the given type.
\end{thm}

\noindent The proof will be given in section \ref{proofs}.

\section{Hecke correspondences} \label{Hecke}

We show how the previous description looks for isogenies of elliptic curves, the case $n=1$.

Recall that on an elliptic curve the polarizations are just the positive integer multiples of the
unique principal polarization, and that an isogeny of elliptic curves acts on polarizations as
multiplication by the degree. With respect to the notation in the previous sections, this means that
if we take $E=(1)$, the principal type, then necessarily $D=(d)$, the degree of the isogeny, and
condition (\ref{tipi}) requires that $\det{M} = \pm d$. Recall moreover that the symplectic group is
the full $SL(2,\Z)$, independent of the polarization degree, and that every $M$ is reduced
under unimodular transformations to a unique diagonal form with diagonal $(d_1,d_2)$ consisting of
positive integers such that $d_1 | d_2$, and necessarily $d_1d_2=d$. With this choice of data we
simply write $\IS(d_1,d_2)$ instead of $\IS_\delta[M]$ for the moduli space of isogenies of the
given type.

This implies for the kernel of the isogeny, occurring in (\ref{seq}), an isomorphism $F\cong
\Z_{d_1}\times\Z_{d_2}$. If $d_1=1$ we have a cyclic subgoup of the elliptic curve which
dominates in the isogeny, and so we see that the space $\IS(1,p)$ is identified to the modular curve
$X_0(p)$ parametrizing elliptic curves with distinguished cyclic subgroup of order $p$. If $d_1>1$
then the isogenies factor as $X \rightarrow Z \rightarrow Y$, an isogeny of type $(1,d_1)$ followed
by an isogeny of type $(1,d_2/d_1)$. 

Let $a,b,p$ be positive integers such that $a|b$ and $(b,p)=1$, and define $d_1=a$ and
$d_2=bp$. Because of the isomorphism $\Z_{d_1} \times \Z_{d_2}  \cong (\Z_a
\times \Z_b) \times \Z_p$, every isogeny $f$ of type $(d_1,d_2)$ factors as $X
\stackrel{u}{\rightarrow} U \stackrel{g}{\rightarrow} Y$, where $u$ is of type $(a,b)$ and $g$ is of
type $(1,p)$, and also factors as $X \stackrel{h}{\rightarrow} V \stackrel{v}{\rightarrow} Y$ where
$h$ is of type $(1,p)$ and $v$ is of type $(a,b)$. In this way two morphisms $\IS(d_1,d_2)
{\rightarrow} X_0(p)$ are defined, sending $f\mapsto g$ and $f\mapsto h$.
The diagram 
$$\begin{array}{ccc} & \IS(d_1,d_2) & \medskip \\ 
& \swarrow\hspace{10pt}\searrow & \\
X_0(p) & & X_0(p)   \end{array}$$
is known as a Hecke correspondence on the modular curve of degree $p$.  This
correspondence is studied in \cite{KM}, \cite{Mi}.
We end with recalling two more examples which fall within the scope of the present treatment.

\begin{exe} \em
{Level structures.}  If $X$ is a variety with polarization
of type $D$, there is an isogeny $\phi:X\rightarrow\hat{X}$ onto the dual variety. A structure on $X$ of
level $D$ is a symplectic basis of \ $\Ker(\phi)$. This is represented by a
symplectic basis $\lambda$ of $H_1(X)$ such that 
$\phi_1(\lambda) (D^{-1} \times D^{-1})$ is a basis of $H_1(\hat{X})$.
The basis of $H_1(\hat{X})$ is symplectic for a principal polarization, type $E=1$. 
With respect to the bases the isogeny $\phi$ is represented by the matrix $M=D\times D$. 
In other words the space of polarized varieties with
structure of level $D$ is a discrete quotient of the space ${\ISP}_{(D,1)}(D\times D)$. 
The quotient is by the discrete group which produces in $X$
a change of symplectic basis without changing the symplectic basis in $\Ker(\phi)$.
See \cite[Ch. 8, \S 3]{LB}.
\end{exe}

\begin{exe} \em
{Varieties with isogeny.}  If $X$ is a variety with
polarization of type $D$, an isogeny $f:X\rightarrow Y$ such that $\Ker(f)$ admits as 
elementary divisors the elements of the diagonal $D$ is called an isogeny of type $D$. 
A canonical basis of $\Ker(f)$
is represented by a symplectic basis $\lambda$ of $H_1(X)$ such that $f_*(\lambda) (1\times D^{-1})$
is a basis of $H_1(Y)$. The basis of $H_1(Y)$ is a symplectic basis
for a principal polarization on $Y$, from which the polarization on $X$ is obtained via pullback. 
With respect to the bases
the isogeny $f$ is represented by the matrix $M=1\times D$. In other words the space of abelian 
varieties  with isogeny of type $D$ coincides with the space ${\IS}_{(D,1)}[1\times D]$. 
See \cite[p. 245]{LB}. The modular variety of abelian surfaces with isogeny of type $D = (1,p)$ is
studied in \cite{HW}, \cite{BL}. 
\end{exe}

\section{Refined invariant from reducibility} \label{refinedinvariant}

For embeddings of abelian varieties a refined invariant is derived from the well known \bigskip

\noindent{\bf Poincar\'e reducibility theorem.} 
{Let $e:X\hookrightarrow Y$ be an embedding into a polarized abelian variety. 
There is an abelian subvariety 
$X'\subset Y$ such that the sum $s:X\times X' {\rightarrow} Y$
is an isogeny. In other words $X\cap X'$ is a finite subgroup and $X+X'=Y$.
Moreover if $L_X := e^* L_Y$ and $L_{X'} := {e'}^* L_Y$ are the induced polarizations, then
$s^* L_Y = p^* L_X \otimes p'{}^* L_{X'}$ in the Neron-Severi group of $X\times X'$, where $p,p'$
are the two projections. Because of this property  the complementary variety $X'$ is
uniquely determined.}
\bigskip

\noindent {\bf Definition.}
We say that an isogeny $s:X\times X' \rightarrow Y$ is a sum of embeddings
if so are the restrictions $X,X' \rightarrow Y$. 
If $Y$ is polarized we say that the isogeny is a
{\em sum of complementary embeddings} if moreover 
$s^* L_Y = p^* L_X \otimes p'{}^* L_{X'}$ holds, where
$L_X, L_{X'}$  are the induced polarizations. 
\bigskip

The reducibility theorem establishes a 1-1 correspondence between embeddings
$e:X \rightarrow Y$ into a polarized abelian variety and isogenies $s:X\times X' \rightarrow Y$
which are sum of complementary embeddings, and this allows to describe the space of embeddings as
the space of isogenies of the special form.
Isomorphisms  of embeddings $e\cong e'$
will correspond with the natural definition of  isomorphisms
$s\cong s'$ of isogenies of the special form.

Let $D,D'$ and $E$ be the types of the polarizations on $X,X'$ and $Y$.
Introduce symplectic bases in $X,Y$ and also in $X'$. 
The induced homomorphism \ 
$s_*:H_1(X\times X') \longrightarrow  H_1(Y)$ is represented by an exact sequence
\begin{eqnarray} \label{seq2}
0 \longrightarrow \Z^{2n}\times\Z^{2n'} \stackrel{M}\longrightarrow
\Z^{2(n+n')} \longrightarrow F \longrightarrow 0
\end{eqnarray}
where $M$ has nonzero determinant, and where $F$ is defined as the cokernel.
Note that the invariant of section \ref{firstinvariant} is the restriction   
$\Z^{2n} \rightarrow \Z^{2(n+n')}$.

The relation $p^* L_X \otimes p'{}^* L_{X'} = s^* L_Y$ of the polarizations is equivalent to
\begin{equation} \label{tipii}
\widehat{D}\times\widehat{D}' = \tp M \widehat{E} M
\end{equation}
where the left hand side is a diagonal block matrix.

We call the datum $\delta=(D,D',E)$
the {\em polarization type} of the embedding, and we call the full datum $\delta, M$ the {\em type of
the embedding relative to the symplectic bases}. 
\bigskip

\noindent {\bf Definition.}
Let $\delta = (D,D',E)$ be a sequence of polarization types, of dimensions $n,n'$ and $n+n'$.
Define $\EM'_{\delta}$ to be the collection of embeddings associated to isogenies 
$s$ which are  sum of complementary embeddings, of the given polarization
type, all varieties endowed with symplectic bases, modulo isomorphisms $s \cong s'$ which preserve
polarizations and bases. 
Define $\EM_{\delta}$ to be the collection of embeddings associated to
isogenies  which are  sum of complementary embeddings,  of the given polarization type,
modulo isomorphisms (which preserve the polarizations). Define
$\Gamma_{\delta} = Sp(D,\Z)\times Sp(D',\Z)\times Sp(E,\Z)$.
\bigskip

The group $\Gamma_{\delta}$ acts on $\EM'_{\delta}$ by producing the change of 
symplectic bases. The element $(A,A',B)$  acts sending the isomorphism class represented
by the isogeny $s$ endowed with symplectic bases $\gamma, \gamma', \lambda$ into the isomorphism
class represented by $s$ with the bases $\gamma A, \gamma' A'$ and $\lambda B$.
There is an identification
$${\EM}_{\delta} \ = \ \begin{array}{c}{\EM}'_{\delta} \\ \hline \Gamma_{\delta}
\end{array}$$
We have seen that there is a natural map
$$\EM'_{\delta} \ \longrightarrow \ GL(2(n+n'), \Z)$$
which to the isomorphism class of $(s, \gamma, \gamma', \lambda)$ associates the matrix
$M$. 

The group
$\Gamma_{\delta}$ also acts on the set of matrices, by producing the equivalence of matrices for a
given homomorphism
$s_*$ under a change of symplectic bases. The element
$(A,A',B)$  acts sending $M \mapsto BM(A\times A')^{-1}$.
The classifying map above is equivariant with respect to the actions of $\Gamma_{\delta}$ and therefore  
induces a map
$$\EM_{\delta} \ = \ \begin{array}{c}\EM'_{\delta} \\ \hline \Gamma_{\delta}\end{array}
\ \longrightarrow \  
\begin{array}{c} GL(2(n+n'), \Z)
\\ \hline \Gamma_{\delta} \end{array}$$
which to the isomorphism class of $s$ associates the equivalence class $[M]$. 
The datum $\delta, [M]$ is what we call the {\em type of the embedding}.

 The fibre over an equivalence class $[M]$ is  described as 
$$\EM_{\delta}[M] \ = \ \begin{array}{c} \EM'_{\delta}(M) \\ 
\hline \Gamma_{\delta}(M) \end{array}$$
where $\Gamma_{\delta}(M)$ denotes the stabilizer of $M$.

If the fibre $\EM'_\delta(M)$ is nonempty the datum $\delta, M$ is an effective {type of embeddings}.
In addition to equality (\ref{tipii}), it is clear that the matrix
$M$ has to satisfy one more property, the property that any isogeny represented by $M$
will be a sum of two embeddings. 
We have the following characterization. 

\begin{lem} \label{property} 
The additional condition so that the datum $\delta, M$ is a type of embeddings is that $M$ fits into
a diagram
\begin{equation*} \label{diagra} \begin{CD} 
0 @>>> \Z^{2n}\times\Z^{2n'} @>{R\times R'}>> \Z^{2n}\times\Z^{2n'}
@>>> F\times F @>>> 0 \\ 
&& || && @AA{}A  @AAA \\  0 @>>> \Z^{2n}\times\Z^{2n'} @>>{M}> \Z^{2(n+n')} @>>> F @>>> 0 
\end{CD} \end{equation*}
where the bottom line is the pullback of the top line under the diagonal homomorphism of $F$, and
where the top line is a product of two exact sequences. 
\end{lem}

\begin{proof}
The isogeny is a sum of embeddings if and only if the inclusion $F\hookrightarrow
X\times X'$ is given by a pair of inclusions $F\hookrightarrow X,X'$.
If the condition is satisfied, as in the beginning of the section, 
calling $\bar X=X/F, \;\bar X'=X'/F$ the quotient varieties, 
the product isogeny $X\times X'\rightarrow \bar X\times\bar X'$ factors as 
$$\begin{CD} X\times X' @>>> \bar X\times\bar X'Ê\\
{||} &&  @AAA \\ X\times X' @>>> Y\end{CD}$$
Introducing bases also in the homology of the
varieties $\bar X$ and $\bar X'$, the diagram is obtained.
Conversely if a diagram of integer homomorphisms as above exists for $M$
then the inclusion of the subgroup into the product variety is given by a pair of inclusions.
\end{proof}

\section{Moduli of embeddings} \label{embeddings}

The refined invariant introduced in the preceding section is sufficient for classifying the
embeddings.

\begin{prop} \label{inclusioni}
For every  type of embeddings there is a natural bijection
$$\EM'_{\delta}(M) \leftarrow\!\rightarrow \AV'_{D}\times\AV'_{D'}$$
\end{prop}

\begin{proof}
To an embedding is associated a pair of varieties, as we have seen.
Conversely, given two abelian varieties $X,X'$, with symplectic bases for polarizations of 
types $D,D'$, by means of the isomorphism
$H_1(X\times X') \cong \Z^{2n}\times\Z^{2n'}$, from the exact sequence (\ref{seq2}) 
an inclusion $F\subset X\times X'$ is deduced. 
Define the torus $Y=X\times X'\,/F$ and call $s:X\times X'\rightarrow Y$ the quotient isogeny.
Then one has  an isomorphism $H_1(Y) \cong \Z^{2(n+n')}$, 
a symplectic basis for an alternating form $L_Y$ of  type $E$. 
Condition (\ref{tipii}) on the types implies that 
$p^*L_X \otimes p'{}^*L_{X'} = s^* L_Y$. Using the same argument as in the proof of Proposition
\ref{duebiiettive} it is then seen that $s^* L_Y$ being a polarization implies that $L_Y$ is a
polarization. Finally because of Lemma \ref{property} the isogeny $s$ is the sum of two embeddings.
In particular  the restriction $X\rightarrow Y$ is an embedding of the given type.
\end{proof}

The stabilizer $\Gamma_\delta(M)$ consists of the triplets 
$(A, A',B) \in \Gamma_\delta$ such that $M(A\times A')=BM$. 
In other words there is a diagram of inclusions
$$\begin{array}{ccc}
&\Gamma_{\delta}(M) \smallskip \\ & \!\!\!^\alpha\swarrow \hspace{20pt} \searrow \\
GL(2n,\Q)\times GL(2n',\Q) & 
\begin{array}{c}{-\!\!\!\longrightarrow}\vspace{-10pt}\\ _{\phi_M}\end{array}
& GL(2n+2n',\Q)
\end{array}$$
where the descending arrows are the natural projections, and
the horizontal arrow is given by $(A,A') \mapsto M (A \times A') M^{-1}$.
Equality (\ref{tipii}) implies that $Sp(D,\Q)\times Sp(D',\Q) =
\phi_M^{-1}\,Sp(E,\Q)$. Therefore
$$\alpha\,\Gamma_{\delta}(M) \ = \ Sp(D,\Z)\times Sp(D',\Z)\;\,\cap\;\,\phi_M^{-1}\,GL(2n+2n',\Z)$$

A product of Siegel spaces $\SH_n\times\SH_{n'}$ is the moduli space of embeddings
with symplectic bases of polarizations of the given types. The group $\alpha\;\Gamma_{\delta}(M)$ 
acts on the product through the actions of symplectic groups. 
Hence the action is properly discontinuous and  the quotient variety exists.
The bijection 
$$\EM_{\delta}[M] \leftarrow\!\rightarrow 
\begin{array}{c}\SH_n\times\SH_{n'} \\ \hline \alpha\,\Gamma_{\delta}(M) \end{array}$$  
introduces on the left hand side the structure of a complex analytic variety. It is an irreducible
variety.

\begin{thm}\label{immerso}
The variety $\EM_{\delta}[M]$ is a coarse moduli space for  embeddings 
of the given type. 
\end{thm}

\noindent The proof will be given in section \ref{proofs}.

\begin{exe}\label{exe} \em
Embeddings of elliptic curves into principally polarized abelian surfaces, the case $n= n'=1$
and $E=I_2$. \hyphenation{a-be-lian} In this situation the pullback polarization
is of type $D=(k)$ where $k=X\cdot\Theta_Y$ is the degree relative to the polarization. Similarly one
has $D'=(k')$ and necessarily $\det M = \pm kk'$. This is seen for instance in
\cite{Mu}, where one also finds some first step in the classification of the  types of embeddings
(Lemma 1) to the effect that the matrix $M$ may be reduced to a form
$$\left( \begin{array}{crcc}
0 &\!\! -k & \cdot & \cdot \\ 0 & 0 & \cdot & \cdot \\ 
1 & 0 & \cdot & \cdot \\ 0 & 1 & \cdot & \cdot  \end{array} \right)$$

The moduli space $\EM_{\delta}(M)$ is a surface, a discrete quotient of
$\SH_1\times\SH_1$. 
Consider the natural map $\EM_\delta(M) \rightarrow \AV_2$. The space of principally
polarized abelian surfaces which contain some elliptic curve of degree $k$ is known to be
irreducible, and is known to coincide with the so called Humbert modular surface of invariant
$\Delta = k^2$. See for instance
\cite{K}, which also contains the historical references. This means that for fixed $k$ all
$\EM_\delta(M)$ have the same image in $\AV_2$.
\end{exe}

\section{Poincar\'e decomposition of morphisms} \label{reducibility}

The treatment of general morphisms is done by patching together the arguments developed so far
for isogenies and embeddings. The starting point is a {Poincar\'e decomposition
of morphisms}.

\begin{prop}
Let $f:V\rightarrow W$ be a morphism of abelian varieties, endowed with polarizations. 
There are complementary abelian subvarieties $X,X'\subset V$  and 
$Y,Y'\subset W$, and  
an isogeny $g:X\rightarrow Y$ such that the following is a commutative diagram
\begin{equation} \label{decomposizione} \begin{CD} 
X\times X' @>>> V \\ @V{g\times 0}VV @VV{f}V \\ Y\times Y' @>>> W
\end{CD} \end{equation}
where the horizontal arrows are isogenies, sums of the given embeddings. 
The diagram is uniquely determined from the polarizations.
\end{prop}

\begin{proof}
Take $X' = \Ker_0(f)$, the connected component of 0 in the kernel, take $Y = f(V)$, let $X$ and
$Y'$ be the complementary abelian subvarieties, and let $g$
be the restriction of $f$.  Uniqueness of the choice is a consequence of
uniqueness in the classical reducibility theorem.
\end{proof}

\noindent{\bf Definition.} In the decomposition above 
if $X'\neq 0$ it is not possible to have on $V$ a pullback polarization. 
However it makes sense to point out the case
in which the isogeny $g:X\rightarrow Y$ preserves the polarizations.
In this situation we say that $f$
is {\em compatible with the polarizations}, or simply a {\em morphism of polarized abelian
varieties}. \bigskip

The proposition establishes a 1-1 correspondence between morphisms $f:V \rightarrow W$ which are
compatible with the polarizations, and sequences of the form
$$\begin{CD} 
X\times X' @>>> V \\ @V{g\times 0}VV  \\ Y\times Y' @>>> W
\end{CD}$$
where the horizontal arrows are isogenies, each a sum of complementary embeddings,
and $g$ is an isogeny which preserves the polarizations, and moreover the sequence is of some
special type, such that a morphism $f$ actually exists which fills in a commutative
diagram  (\ref{decomposizione}). This allows to describe the space of morphisms as the space of
sequences of the special type.

Let $D,D',E$ and $H,H',K$ be the polarization types of  $X,X',V$ and $Y,Y',W$.  
Introducing symplectic bases in all these varieties,
the associated diagram of homology groups is represented by a diagram 
\begin{equation} \label{diagrammo} \begin{CD}
0 @>>> \Z^{2n}\times\Z^{2n'}  @>{M}>> \Z^{2(n+n')} @>>> F @>>> 0
\\ && @V{P\times 0}VV  @VVV @VVV  \\ 
0 @>>> \Z^{2m}\times\Z^{2m'} @>>{N}> \Z^{2(m+m')} @>>> G @>>> 0 \\
\end{CD} \end{equation}
where $n=m$ and where the  matrices $M,N$ and $P$ have nonzero determinant.

The cokernels $F,G$ are isomorphic to the kernels of the horizontal isogenies in the decomposition,
i.e. the subgroups $X\cap X'$ and $Y\cap Y'$. Thus the morphism $f$ is isomorphic to a morphism
$$\begin{array}{c} X\times X' \\ \hline F\end{array} \ \ {\longrightarrow} \ \
\begin{array}{c} Y\times Y' \\ \hline G\end{array}$$
which is a quotient of the product  $g\times 0$  in the diagram.

We denote by
$\delta \ = \ (D,D',E,H,H',K)$
the sequence of types of polarizations, and call it the {\em polarization type} of the morphism.  
We denote by
$\tau \ = \ (M,N,P)$ the sequence of matrices which represent the
decomposition of the morphism, and we call the datum $\delta, \tau$ the {\em type of the morphism
relative to the symplectic bases}. 

The types $D,D',E$ are related to $M$ in equality (\ref{tipii}), 
and similarly the types $H,H',K$ are related to $N$.
As the morphism is compatible with the polarizations then
$\widehat{D} =\tp P\widehat{H}P$ also holds as in (\ref{tipi}).
The matrices $M,N,P$ satisfy a number of conditions. 
The product $N\,(0\times P)M^{-1}$ is an integer matrix, associated to $f$. 
Both matrices $M,N$ satisfy the condition of Lemma \ref{property}. 
Finally it is quite clear that one more property holds. It is the property that, if a sequence of
isogenies of the special form is given, and if the sequence is of the given type $\delta,\tau$, then
a morphism $f$ actually exists which fills in the commutative diagram (\ref{decomposizione}). 
A characterization is given in the following lemma. If all these properties are verified then 
the datum $\delta,\tau$ is an effective  {type of morphisms}.

\begin{lem} \label{strongproperty} 
The additional condition for the datum $\delta,\tau$ as above to be a type of morphisms is
that the matrix $P$ fits in some diagram
\begin{equation*} \label{diagra} \begin{CD} 
0 @>>> \Z^{2n} @>{R}>> \Z^{2n} @>>> F @>>> 0 \\   && || && @VV{\bar P}V   \\  
&& \Z^{2n} @>>{P}> \Z^{2n}   \end{CD} \end{equation*}
in other words that $P = \bar P R$ with  $\mbox{\rm Coker}(R) = F$.
\end{lem}

\begin{proof}
By hypothesis the matrices of the isogenies satisfy the condition of Lemma 
\ref{property} and there are inclusions $F\subset X,X'$ and $G\subset Y,Y'$ 
as is seen in the proof of the lemma. 
Using these inclusions it is easily seen that there exists $f:V\rightarrow W$ 
if and only if $g\times 0$ sends $F\rightarrow G$ and this happens if and only if 
$g:X\rightarrow Y$ sends $F\rightarrow 0$.
If this happens then there is a factorization $X\rightarrow\bar X\rightarrow Y$
where $\bar X =X/F$. Introducing bases in homology as in the proof of the lemma,
the associated matrices satisfy $P=\bar P\,R$ where $\bar P$ 
is the matrix of $\bar X\rightarrow Y$ and where $R$ is the matrix of 
$X\rightarrow\bar X$. 
Conversely if the matrix $\bar P$ exists then
$F\rightarrow 0$ under $X\rightarrow Y$.
\end{proof}

\section{General moduli spaces} \label{generalmorphisms}

Finally we show that the moduli spaces of morphisms of  given types are irreducible varieties. 
\bigskip

\noindent {\bf Definition.}
Let $\delta$ be a sequence of 6 polarization types, whose dimensions are in the relation as
in the preceding section, diagram (\ref{diagrammo}). Define $\AVAP_{\delta}$ to be the collection of
morphisms
$f$ admitting a decomposition (\ref{decomposizione})
of the given polarization type, all varieties endowed with symplectic bases,
modulo isomorphisms which preserve polarizations and bases. Define 
$\AVA_{\delta}$ to be the collection of morphisms
admitting a decomposition of the given polarization type,
modulo isomorphisms (which preserve the polarizations). Define 
$\Gamma_\delta$  to be the product of the 6 symplectic groups 
for the types in $\delta$.
\bigskip

The group $\Gamma_\delta$ acts on $\AVAP_{\delta}$ by producing the change of symplectic bases,
and there is an identification
$$\begin{array}{ccc} 
{\AVA}_{\delta} &=& \begin{array}{c} \AVAP_{\delta} \\ 
\hline \Gamma_{\delta} \end{array} \end{array}$$
For the natural map on $\AVAP_{\delta}$ which to the isomorphism class of a morphism $f$
endowed with bases associates the type $\tau$ relative to the bases, we denote by 
$\AVAP_{\delta}(\tau)$ the fibre over $\tau$. It is nonempty if the datum $\delta,\tau$ is a type
of morphisms.

The group $\Gamma_\delta$ also acts on the set of sequences $\tau$ by producing the equivalence 
under the change of bases. The classifying map on $\AVAP_\delta$ is equivariant, and therefore
defines a map on $\AVA_\delta$ which to the isomorphism class of $f$ associates the equivalence
class $[\tau]$. The datum $\delta, [\tau]$ is what we call the {\em type of the morphism}. 

If $\AVA_\delta[\tau]$ is the fibre over $[\tau]$, we have the natural identification

$$\begin{array}{ccc} 
{\AVA}_{\delta}[\tau] &=& 
\begin{array}{c} \AVAP_{\delta}(\tau) \\ 
\hline \Gamma_{\delta}(\tau)\end{array} \end{array}$$
where  $\Gamma_{\delta}(\tau)$  is the subgroup of 
$\Gamma_\delta$ stabilizer of $\tau$.

\begin{prop} \label{sopra}
For every type of morphisms
there is a natural bijection
$$\AVAP_{\delta}(\tau) \ \longleftrightarrow \ 
\ISP_{D,H}(P)\times \AV'_{D'} \times \AV'_{H'}$$
\end{prop}

\begin{proof}
Assume that we are given abelian varieties $X,X',Y,Y'$ with symplectic bases for polarizations of 
 types $D,D',H,H'$, and an isogeny $g:X\rightarrow Y$ of  type $(D,H),P$.
Because of Lemma \ref{property} for $M$ there are inclusions $F\subset X,X'$. Then define the
torus $V=X\times X'\,/F$. By means of diagram (\ref{diagrammo}) an isomorphism 
$H_1(V) \cong \Z^{2(n+n')}$ is obtained. Because the types $D,D',E$ are related, this gives on
$V$ a basis for a polarization of  type $E$, as follows from Proposition \ref{inclusioni}.
Similarly there are inclusions $G\subset Y,Y'$ and defining $W=Y\times Y'\,/G$ one has an
isomorphism  $H_1(W) \cong \Z^{2(m+m')}$, a basis for a polarization of  type $K$. 
From Lemma  \ref{strongproperty} 
it follows that the morphism $g\times 0:X\times X'\rightarrow Y\times Y'$ sends $F\rightarrow G$
and  induces a morphism $f:V\rightarrow W$. With respect to the bases
this morphism is of the given type $\delta,\tau$.
\end{proof}

It was seen in \ref{locclos} that the family $\ISP_{D,H}(P)$ is represented by the Siegel space
$\SH_{n}$. It follows from the proposition above that  $\AVAP_\delta(\tau)$ is
represented by the product $\SH_n \times \SH_{n'} \times \SH_{m'}$, and therefore there is a
bijective map
$$\AVA_\delta[\tau] \longleftrightarrow \frac{\SH_n \times \SH_{n'} \times \SH_{m'}}
{\Gamma_\delta(\tau)}$$ 
It is seen in the usual way that the   stabilizer
$\Gamma_\delta(\tau)$ acts properly discontinuously, so the quotient variety exists, and this
introduces on the left hand side the structure of a complex analytic variety, and an
irreducible one.

\begin{thm}\label{generale}
The variety ${\AVA}_{\delta}[\tau]$ is a coarse moduli space for morphisms of the given type.
\end{thm}

\noindent The proof will be given in section \ref{proofs}.

\section{The moduli property}\label{proofs}

In this section we give the proof of the moduli
property, in a sketched form. The proof is
based on an analytic treatment of the moduli
property of the Siegel modular variety. 

The moduli property means two things. First,
each modular variety parametrizes 
\hyphenation{pa-ra-me-tri-zes}
a family of morphisms. 
This is constructed in a natural way from
the universal family of abelian varieties 
over the Siegel space. The details are
lenghty, but natural, and will be omitted.
Second, every flat family of morphisms 
determines a morphism to the suitable variety 
of moduli. This part of the proof will be
sketched in the various cases.

\bigskip
\noindent{\it Proof of } \ref{isogeniale}.
A flat family of isogenies is a surjective morphism
$f:\X \rightarrow \Y$
between two flat families of polarized 
abelian varieties of the same dimension over
some scheme $T$.  Locally over open subschemes 
$T'\subset T$ the two families can be endowed
with symplectic bases and this determines two
liftings
$T'\rightarrow \AV'_D,\, \AV'_E$. With
respect to the bases all isogenies in the
family are represented by some constant integer
matrix, which we may assume is the given $M$.
Then the two liftings above are two morphisms
$T'\rightarrow \SH_n$ which are
identified under the automorphism  
$\SH_n \leftrightarrow \SH_n$ of Proposition \ref{locclos},
and so both represent on $T'$ the same  
morphism  $T\rightarrow {\IS}_{\delta}[M]$.
\qed \bigskip

\noindent{\it Proof of } \ref{immerso}.
A flat family of embeddings
$\X\rightarrow\Y$ into a polarized
family of abelian varieties determines a flat
family of complementary varieties $\X'$, with
the induced polarizations. Locally over
open subschemes $T'$ of the parameter scheme
$T$ the families can be endowed with symplectic
bases, and this determines two liftings
$T'\rightarrow\AV'_D,\AV'_{D'}$ and so a
morphism 
$T'\rightarrow\SH_n\times\SH_{n'}$ which locally represents a morphism 
$T\rightarrow\EM_{\delta}[M]$. 
\qed \bigskip

\noindent{\it Proof of } \ref{generale}.
Let $f:\V\rightarrow\W$ be a flat family of morphisms of the given type, 
\hyphenation{pa-ra-me-tri-zed} 
parametrized by $T$. By hypothesis all
morphisms in the family are of the same rank.
Because of this, there is a flat subfamily
$\X'\subset\V$ which parametrizes the 
connected kernels $\Ker_0(f_t)$, and the image 
$\Y=f(\V)$ also is a flat family. Then the
families $\X,\Y'$ of complementary varieties
are determined relative to the polarizations,
and finally is defined a family of isogenies
$g:\X\rightarrow\Y$. They are all flat
families.
Locally over open subschemes
$T'\subset T$ the families of abelian
varieties can be endowed with symplectic bases
for the type $\delta$, with respect to which
the morphisms occurring in the  factorization
of $f_t$ are represented by a sequence of
integer matrices, independent of $t$,
belonging to the given type, 
we may assume it is precisely the sequence
$\tau$. From the universal property of the
Siegel space and from Theorem 1.1 it follows
that there is a morphism
$T'\rightarrow \ISP_{D,H}(P)\times \AV'_{D'}\times \AV'_{H'}$
which locally represents a morphism
$T\rightarrow\AVA_\delta[\tau]$.
\qed \bigskip

Finally we may define the total space $\AVA  = \, \cup  \ \AVA_\delta[\tau]$,
disjoint sum of countably many analytic varieties.
This will be a coarse moduli space for all morphisms
of polarized abelian varieties.
What is only required in addition is the
observation that in a family of morphisms
$f:\V\rightarrow\W$ the image $f(\V)$ has constant fibre dimension over
the parameter $T$. This is a general statement of rigidity for morphisms of abelian varieties,
and will be considered elsewhere. 
Then as in the preceding
proof the decomposition is constructed,  with
respect to which the family is of some well
determined type
$\delta,[\tau]$,  and finally a morphism
$T\rightarrow\AVA_\delta[\tau]$ is obtained.

\section{Some questions in integral linear algebra} \label{questi}

We have seen that the components of the moduli space of morphisms between polarized abelian varieties
correspond to  equivalence classes of certain rather complicated discrete data. Thus the question
arises of finding canonical forms for these data. All what we know  is contained in
the examples of section \ref{Hecke} and in example \ref{exe}. Here we recollect some results and
formulate some natural questions.

The types of isogenies are the data of the form $\delta, M$ where $\delta = (D,E)$ is a pair
of polarization types of the same dimension $n$, and $M$ is an integer square matrix of order $2n$,
such that  $\widehat{D}  =  \tp M \widehat{E}M$ holds (equality (\ref{tipi})). 

For given $\delta$ we ask whether some $M$ exists. A necessary condition is that $E$ divides $D$, as
follows from \cite[Ch. 9, p. 85, ex. 1a]{B}. Moreover, a solution $M$ is equivalent to  $BMA^{-1}$
for every element
$(A,B)$  of the group $\Gamma_\delta$,  a pair of symplectic matrices for the types
$D,E$. It would be useful to find canonical forms, i.e. the simplest possible forms in the
equivalence classes $[M]$. Finally, it is clear that equality (\ref{tipi}) says that the smaller
datum $E,M$ determines the datum $D$. It might be useful to rephrase the whole treatment in
terms of the smaller datum.

The types of embeddings are the data of the form $\delta, M$ where $\delta = (D,D',E)$ is a sequence
of polarization types, of dimensions $n,n',n+n'$, and $M$ is an integer square matrix of order
$2(n+n')$, which is obtained from some product $R \times R'$ in the way described in Lemma
\ref{property}, and such that 
$\widehat{D} \times \widehat{D'}  = \tp M \widehat{E}M$ holds (equality (\ref{tipii})). 

For the same reason as before we find that a necessary condition is that $E$ divides the
type of the matrix $\widehat{D} \times \widehat{D'}$, which is determined by the pair 
$D,D'$, cf. \cite[Ch. 7, p. 94, rem. 3]{B}. Moreover, we may replace $M$ with any equivalent matrix
of the form $BM(A \times A')^{-1}$, for any element $(A,A',B)$  of the group $\Gamma_\delta$, and it
would be useful to find some simpler equivalent form. It is possible that this problem, for  
embeddings, might become easier to work with once we have some solution of the previous problem, for
isogenies.

The analogous discussion can be done about the classification of the general types of
morphisms.

\vfill
\noindent
address: \ Dipartimento di Matematica e Informatica, Universit\`a di Perugia, 
via Vanvitelli 1, 06123 Perugia, Italy \smallskip\newline
e-mail: \ {\tt guerra@unipg.it}

\end{document}